\documentclass[10pt,fleqn]{article}
\usepackage{algorithmic}

\usepackage{latexsym}
\usepackage{amsmath}
\usepackage{amscd}
\usepackage{amssymb}
\usepackage{amsthm}
\usepackage{amsopn}

\theoremstyle{plain}
\newtheorem{theorem}{Theorem}
\newtheorem{Lemma}[theorem]{Lemma}
\newtheorem{Theorem}[theorem]{Theorem}
\newtheorem{Corolary}[theorem]{Corollary}
\newtheorem{Proposition}[theorem]{Proposition}
\theoremstyle{definition}
\newtheorem{Definition}[theorem]{Definition}
\newtheorem{Example}[theorem]{Example}

\newcounter{myenum}[theorem]
\newenvironment{eDefinition}[1][{}]
{
  \begin{Definition}
  {{#1 \mbox{}  }}
  \begin{enumerate}
  \usecounter{myenum}
  \numberwithin{myenum}{theorem}

}
{\end{enumerate}\end{Definition}}

\newenvironment{eExample}[1][{}]
{
  \begin{Example}
  {{#1 \mbox{}  }}
  \begin{enumerate}
  \usecounter{myenum}
  \numberwithin{myenum}{theorem}

}
{\end{enumerate}\end{Example}}

\def\ZZ{{\mathbb{Z}}}
\def\QQ{{\mathbb{Q}}}
\def\B{{\mathcal{B}}}
\def\G{{\mathcal{G}}}
\def\g{{\mathfrak{g}}}
\def\U{{\mathfrak{U}}}
\def\Grobner{{Gr\"{o}bner }}

\newcommand\LM[2][{}]{\mathbf{LM}_{{#1}}\left({#2}\right)}
\newcommand\LC[2][{}]{\mathbf{LC}_{{#1}}\left({#2}\right)}
\newcommand\LT[2][{}]{\mathbf{LT}_{{#1}}\left({#2}\right)}

\title{Unital \Grobner Bases over Aribtrary Commutative Ground Rings}
\author{Frederick Leitner, Robert Pawloski\footnote{Both authors were supported under a NSF VIGRE grant.}}

\begin{document}

\maketitle 

\begin{abstract}
	Let $R$ be  a commutative ring with unity and a let $A$ be a not necessarily
	commutative $R$-algebra which is free as an  $R$-module.
	If $I$ is an ideal in $A$, one can ask when $A/I$ is also
	free as an $R$-module.   We show that if $A$ has an \emph{admissible
	system} and $I$  has a \emph{unital \Grobner
	basis} then $A/I$ is free as an $R$-module.   We prove a version of
	Buchberger's theorem over $R$ and, as a corollary, we
	obtain a \Grobner basis proof of the Poincare-Birkhoff-Witt Theorem over
	a commutative ground ring.

	\noindent MSC: 16Z05,13P10
\end{abstract}

\section{Introduction.}
	There have been several generalizations of \Grobner basis theory,
	coming in one of two flavors:  noncommutative theory and theories for working
	over special ground rings, e.g.   Euclidean domains, PID's or UFD's.
	For an overview of the theory of such rings see \cite{adams}.
	Over a field, one of the main uses of a  \Grobner basis is to find 
	a basis for the quotient of an algebra $A$ by a (left or two-sided) ideal $I$.
	Let us consider two examples for $k$ a (commutative) field (with unity)
	whose quotient algebras can be readily described through the known \Grobner
	basis theory.  
	\begin{eExample}
	\item
		Let $A$ be the polynomial algebra $k[x_1,\cdots,x_n]$ and consider the 
		ideal $m=(x_1,\cdots,x_n)$.  
		Taking
		a graded lexicographic ordering on $A$ with $x_i < x_j$ if $i < j$,
		$m^2$ has a \Grobner basis given by $\{x_i x_j\}$.
		Thus the quotient	$A/m^2$ is an $n+1$-dimensional 
		$k$-vector space with a basis given by:
		\begin{displaymath}
			A/m^2=k\{1,x_1,\cdots,x_n \}
		\end{displaymath}

	\item
		If $\g$ is a 
		lie algebra with lie bracket $[ \ , \ ]_{\g}$  over $k$ then one forms the universal
		enveloping algebra $\U\g$ as the quotient of the tensor or free algebra $T\g$
		on $\g$ by the two sided ideal $J$:
		\begin{displaymath}
			J=\left( xy-yx -[x,y]_{\g} \ | x,y \in \g \hookrightarrow T\g \right)
		\end{displaymath}
		If one chooses a total ordering $<$ on the index set $I$ for 
		a $k$-basis $\g=k<x_i | i \in I>$ then a \Grobner basis argument
		yields the Poincare-Birkhoff-Witt(PBW) theorem which says that $\U\g$ has a basis
		of non-decreasing words (\cite{Mora} \cite{DeGraaf}).
	\end{eExample}
		
	\noindent A consideration of the first example shows the fact that $k$ was a field
	was not important --  one would have a similar statement with the
	integers $\ZZ$, though not through \Grobner basis techniques.   The same holds for
	the second example  if we replace $k$ by an arbitrary (commutative) ring.
	However, in this more general setting, one cannot make use of \Grobner basis techniques
	-- one must prove this by ``ad-hoc''  methods as in \cite{Serre}.  
	However, these arguments  bear a close relation to those used in the \Grobner basis 
	theory.  We view the inability of \Grobner basis techniques to apply to these 
	mild generalizations  as an unsatisfactory state of affairs.

	If one regards the above two examples closely, one sees that the ground field $k$
	never enters into the picture.  Specifically, one does not need to invert any constants, 
	and this leads to the notion of a \emph{unital \Grobner basis}.
	We will prove the following:

	\begin{Theorem}\label{thm}
		Let $R$ be a commutative algebra with unity.
		Let $A$ be an $R$-algebra which is  free as an $R$-module and without 
		quasi-zeros.
		Let $(\B,<)$ be an admissible system on $A$ and let $I$ be 
		a two-sided  ideal of $A$ with a unital \Grobner basis $\G$.
		Define $\tilde{O}(\G)$ to be  the free $R$-module spanned
		by the monomials which do not occur as leading monomials  of members of $\G$.
		Then:
		\begin{enumerate}
		\item\label{thm1}
			There is a $k$-module isomorphism $A/I\simeq \tilde{O}(\G)$ of free  $R$-modules.
		\item\label{thm2}
			$A=I\oplus \tilde{O}(\G)$ 
		\end{enumerate}
	\end{Theorem}

	\noindent For the case of a left ideal $I$, or for the case where $A$ has quasi-zeros, one only need
	to combine the techniques in \cite{Li} with ours.
	As a corollary we obtain the PBW theorem over an arbitrary commutative ring with unity.
	In particular, if $R$ is a $\QQ$-algebra we have a \Grobner basis proof of the equivalence 
	of the category of (finite dimensional) lie algebras over $R$ and (finite dimensional) smooth 
	formal groups over $R$.

\section{Unital \Grobner Bases.}
	Throughout this section
	$R$ is a commutative ring with unity and $A$ an $R$-algebra with no quasi-zeros, i.e.
	elements $a\in A$ such that for all $b,c \in A$ not
	both $1$ one has $bac=0$.
	For brevity, we only state and prove the results
	in the case of a two-sided ideal.   
	We take time to fix notation, following closely that of (\cite{Li}), 
	
	\begin{eDefinition}
		{Let $A$ be an $R$-algebra with multiplication $\cdot$ (and without quasi-zeros).
		Choose a set of (algebra) generators
		$A=R\langle x_i | i\in \Lambda\rangle$ for some index set $\Lambda$.}
		\item
			Let $\alpha$ be a finite length word in $\Lambda$, i.e. an ordered expression:
			\begin{displaymath}
				\alpha=\alpha_1 \alpha_2 \cdots \ \alpha_n \ \ \ \alpha_i \in \Lambda
			\end{displaymath}
			Then a  \emph{monomial} in $A$ is the ordered product:
			\begin{displaymath}	
				x^{\alpha}:= x_{\alpha_1} \cdot x_{\alpha_2} \cdots \ x_{\alpha_n}
			\end{displaymath}
		\item
			Suppose that $A$ is a free $R$-module.
			We say that $A$ has a \emph{monomial basis $\B$ sub-ordinate to $\Lambda$}
			if  $\B$ is a subset of all words of finite length  in $\Lambda$
			such that $A$ has  an $R$-basis:
			\begin{displaymath}
				A=R\{ x^{\beta} \ | \ \beta \in \B \}
			\end{displaymath}
			We identify $\beta$ with $x^{\beta}$ to simplify notation.
	  	\item
			Assume that $\B$  is well ordered by $<$.
    			Let $f\in A$, then we may write $f$ uniquely as the finite sum:
    			\begin{displaymath}
				f= c_1 b_1 + c_2 b_2 \cdots + c_n b_n \mbox{ for }
			    	c_i \in R^{\ast}:=R\setminus\{0\} 
				\mbox{ and }b_i \in \B.
    			\end{displaymath}
    			such that 
    			$b_1 > b_2 > \cdots > b_n$.
    			Then  the \emph{leading (or head) monomial of $f$ with respect to $\B$ and $<$} is defined as:
    			\begin{displaymath}
      				\LM{f} := \LM[{<}]{f}:= b_1
    			\end{displaymath}
    			i.e. the largest basis element appearing.
    			The \emph{leading (or head) coefficient of $f$ with respect to $\B$ and $<$} is defined as:
    			\begin{displaymath}
      				\LC{f} := \LC[{<}]{f}:= c_1
    			\end{displaymath}
			while the \emph{leading (or head) term of $f$ with respect to $\B$ and $<$} is defined as:
			\begin{displaymath}
      				\LT{f} := \LT[{<}]{f}:= \LC{f}\LM{f}
			\end{displaymath}
		\item
			A \emph{ monomial ordering} for a basis $\B$ of  an algebra 
    			$B$   is a well ordering $<$ on $\B$
    			such that  for $b,b',r,s \ \in \B $ we have:
    			\begin{enumerate}
    			\item
			      	if  $b < b'$  then  $r\cdot b\cdot s < r  \cdot b' \cdot s$
      				whenever:
      				\begin{displaymath}
					r\cdot b \cdot s \mbox{ and } \LM{r \cdot b' \cdot s} \ne 0
     				\end{displaymath}	
    			\item
      				if  $b' =\LM{r \cdot b \cdot s}\ne 0$  with 
      				$r \mbox{ or } s \ne 1$  then  $b < \LM{b'}$
  			\end{enumerate}
  		\item
    			If   $\B$ admits a  monomial ordering $<$, then the pair
    			$(\B,<)$ is an \emph{admissible system}.
		\item
			If $(B,<)$ is an admissible system, and $0\ne f\in A$, then
			we say that $\beta \in \B$ divides $f$ ( and write $\beta|f$) if
			there are $u,v\in \B$ and $\lambda \in R^{\ast}$ such that:
			\begin{displaymath}
				\lambda \LT{u\cdot \beta \cdot v} = \LT{f}
			\end{displaymath}
	\end{eDefinition}

	\noindent We now describe for  a subset $I\subseteq A$, $I$ not necessarily an ideal,
	a ``division algorithm.''   We have put quotations to emphasize
	that this algorithm does not in general have a meaning.  In fact, one
	may view the statement of (Theorem \ref{thm}) as ascribing a meaning to this
	algorithm when we make the assumption of the existence of a \emph{unital \Grobner basis} 
	for $I$.
	A second reason for putting this in quotations is that we give no prescription
	for choosing the elements of $I$ with which to  divide. However, this need not be a hindrance 
	and is in fact a benefit in view of  (Lemma \ref{lem}).
	We need one more set of  definitions at this point:

	\begin{Definition}
		Let $R$ be a commutative ring with unity, $A$ an $R$-algebra
		which is free as an $R$-module and without quasi-zeros, $(\B,<)$ an admissible system
		on $A$, and a subset $I\subseteq A$.  Then define a 
		$R$-submodule $O(I)$ of $A$ to be the $R$-submodule spanned by
		the set:
		\begin{displaymath}
			o(I) := \{ \lambda \beta \ | \ \lambda \in R^{\ast}, \  \beta \in \B, \ \forall h\in I
			\mbox{ we have  }
			\LT{h} \ne \lambda \beta \}
		\end{displaymath}
		We also define the $R$-module $\tilde{O}(I)$ as the $R$-submodule spanned by
		the set:
		\begin{displaymath}
			\tilde{o}(I) := \B \setminus \{ \LM{h} \ | \ h \in I \}
		\end{displaymath}
		Clearly $\tilde{O}(I)$ is a free $R$-module.
	\end{Definition}

	\renewcommand{\algorithmicrequire}{\textbf{Input:}}
	\renewcommand{\algorithmicensure}{\textbf{Output:}}
	\begin{algorithmic}[1]
		\REQUIRE$R$ a commutative ring with unity,
		$A$ an $R$-algebra which is free as an $R$-module and without quasi-zeros,
		$(\B,<)$ an admissible system on $A$,
		a subset $I\subseteq A$, and  $f\in A$.
		\ENSURE $\tilde{f} \in I$ and 
			$r\in O(I)$ the \emph{remainder of $f$ on division by $I$}
		        so that $f=r + \tilde{f}$ 
		\STATE  $i:=0$ 
		\STATE  $f_0 := f$.
		\WHILE {$f_i \ne 0$}
			\STATE $i:=i+1$
			\IF {$\not\exists h \in I$ such that $\LT{h} | f$}
				\STATE $r_i := \LT{f}$
				\STATE $f_i := f_{i-1} - r_i$
			\ELSE
				\STATE Choose some $h_i \in I$ such that $0\ne\LT{h_i} | f_{i-1}$
				\STATE Choose some $\lambda_i \in R^{\ast}$ $u_i,v_i\in \B$ so that:
				\STATE \hspace{1cm} $\lambda_i\LT{u_i\cdot h_i \cdot v_i} = \LT{f_{i-1}}$
				\STATE $f_i := f_{i-1} - \lambda_i u_i\cdot h_i \cdot v_i$
				\STATE $r_i := 0$
			\ENDIF	
		\ENDWHILE
		\STATE $r:= \sum_{i} r_i$
		\STATE $\tilde{f}:=\sum_i \lambda_i u_i \cdot h_i \cdot v_i = f - r$
	\end{algorithmic}
	
	\noindent We note that, because we do not specify how to choose the $h_i$ (nor the $u_i$ and $v_i$) we do
	not in general have a unique output.

	\begin{eDefinition}
		{Let $R$ be a commutative ring with unity, $A$ an $R$-algebra which is free as
		a $R$-module and without quasi-zeros, and let $(\B,<)$ be an admissible system on $A$.}
		\item\label{gbdef}
		   	Let $I$ be a (two-sided)
			ideal in $I$.  Let $I$ have a set of generators:
			\begin{displaymath}
				\G=\{ g_\gamma \ | \ \gamma \in \Gamma \} \mbox{ for some index set } \Gamma
			\end{displaymath}
			Then we say that $\G$ is a \emph{\Grobner basis with respect to $(\B,<)$} if
			for every $h\in I$ we have a representation:
			\begin{displaymath}
			  h = \sum_{k\in K} \lambda_{k} u_{k}\cdot g_{\gamma_k} \cdot v_{k}
			\end{displaymath}	
			for $K$ an index set, $\gamma_k \in \Gamma$, 
			$\lambda_{k}\in R^{\ast}$ and $u_{k},v_{k} \in \B$
			such that $\LM{u_{k}\cdot g_{\gamma_k}\cdot  v_{k}} \le \LM{f}$ whenever
			$u_{k}\cdot g_{\gamma_k}\cdot v_{k}\ne 0$
		\item
		        We call a subset $\G \subseteq A$ \emph{unital} if:
			\begin{enumerate}
			\item
			  For all $\gamma \in \Gamma$ we have:
			  \begin{displaymath}
			    \LC{g_{\gamma}}\in R^{\times}:= \mbox{ units of } R
			  \end{displaymath}
			\item
			  For all $\gamma \in \Gamma$ and for all $\alpha,\beta \in \B$ we have
			  \begin{displaymath}
			    \LC{\alpha\cdot g_{\gamma} \cdot \beta} \in R^{\times} \mbox{ whenever } \alpha\cdot g_{\gamma} \cdot \beta \ne 0
			  \end{displaymath}	
			\end{enumerate}
		\item
			For $f,f'\in A$ we say that an \emph{S-polynomial is constructible about $f$ and $g$}
			if there is some $u,u',v',v' \in \B$, $\lambda,\lambda' \in R^{\ast}$ such
			that:
			\begin{displaymath}
				\lambda \LT{u\cdot \LT{f}\cdot v} =
				\lambda' \LT{u'\cdot \LT{f'}\cdot v'} 
			\end{displaymath}
			In which case we write:
			\begin{displaymath}
				S:=S(f,f'):= \lambda u \cdot f \cdot v - \lambda' u' \cdot f' \cdot v' 
			\end{displaymath}
			We say that $S$ is an \emph{$S$-polynomial about $f$ and $f'$}.  The choices
			of $u,u',v,v'$ are not in general	unique.
	\end{eDefinition}

	\noindent The extra conditions of  a \Grobner basis being unital are not that strong:
	If the ground ring is a field, then a \Grobner basis is automatically a unital \Grobner basis.
	If $\B$ is closed under multiplication, then the third
	condition follows if the second condition holds.  We also note that the construction of a
	$S$-polynomial 	ensures that:
	\begin{displaymath}
	  \LM{S} < \LM{u\cdot \LM{f}\cdot v} = \LM{u'\cdot \LM{f'}\cdot v'}
	\end{displaymath}

	\begin{Lemma}[\cite{MadRein}]\label{lem}
		Let $R$ be a commutative ring with unity, $A$ an $R$-algebra which is free as an $R$-module,
		$(\B,<)$ an admissible system on $A$, and let $I$ be a two sided ideal of $A$ which
		is generated by a unital \Grobner basis:
		\begin{displaymath}
			\G = \{ g_{\gamma} \ | \ \gamma \in \Gamma \}.
		\end{displaymath}   Then in the division algorithm,
		we may choose $h_i \in I$ so that $h_i = g_{\gamma_i}$ for some $\gamma_i \in \Gamma$.
	\end{Lemma}
	\begin{proof}
		Let us set $h:=h_i$, $f:=f_i$ to stop the proliferation of subscripts.
		Then since $\G=\{ g_{\gamma} \ | \ \gamma \in \G\}$ is a \Grobner basis, we may write:
		\begin{equation}\label{eqn}
			h=\sum_{k \in K} \lambda_k u_k \cdot g_{\gamma_k} \cdot v_k 
		\end{equation}
		with $\lambda_k \in R^{\ast}$, $u_k, v_k \in \B$ and $\gamma_k\in \Gamma$.
		Denote:
		\begin{displaymath}
			\alpha := \LM{h}
		\end{displaymath}
		As we are free to choose our representation (Equation \ref{eqn}) of $h$ as we
		wish, we may choose one so that $\alpha$ is minimal with respect to the ordering
		$<$.    
		Denote:
		\begin{displaymath}
			T:= \{ k \in K \ | \ \LM{u_k \cdot g_{\gamma_k} \cdot v_k} = \alpha \} 
		\end{displaymath}
		We can further choose a representation of $h$ so that $|T|$ is minimal.
		If $|T| =  1$ we are done.  Otherwise, let $k_1\ne k_2\in T$, and denote
		$c_{k_i}:= \LC{u_{k_i} \cdot g_{\gamma_{k_i}}\cdot v_{k_i}}$.   By the
		assumption that $\G$ was a unital \Grobner basis, we have that
		$c_{k_i} \in R^{\times}$ so that we may form the  
		$S$-polynomial:
		\begin{displaymath}
			S:=\lambda_{k_2}\frac{c_{k_2}}{c_{k_1}} u_{k_1}\cdot g_{\gamma_{k_1}}\cdot v_{k_1}
			-\lambda_{k_2} u_{k_2}\cdot g_{\gamma_{k_2}}\cdot v_{k_2}
		\end{displaymath}
		Then we have:
		\begin{eqnarray*}
			\lefteqn{h= 
			    \lambda_{k_1} u_{k_1}\cdot g_{\gamma_{k_1}}\cdot u_{k_1}
			    +\lambda_{k_2} u_{k_2}\cdot g_{\gamma_{k_2}}\cdot u_{k_2}
			    + \sum_{k\ne k_1,k_2} \lambda_{k} u_{k}\cdot g_{\gamma_{k}}\cdot u_{k}}
			\\
			&=&
			    \lambda_{k_1} u_{k_1}\cdot g_{\gamma_{k_1}}\cdot u_{k_1}
			    + \left(  \begin{array}{c}
					\lambda_{k_2}\frac{c_{k_2}}{c_{k_1}} u_{k_1}\cdot g_{\gamma_{k_1}}\cdot v_{k_1} \\
					-\lambda_{k_2}\frac{c_{k_2}}{c_{k_1}} u_{k_1}\cdot g_{\gamma_{k_1}}\cdot v_{k_1}
				      \end{array}
				\right)
			    +\lambda_{k_2} u_{k_2}\cdot g_{\gamma_{k_2}}\cdot u_{k_2}
			\\
			&&
			    + \sum_{k\ne k_1,k_2} \lambda_{k} u_{k}\cdot g_{\gamma_{k}}\cdot u_{k}
			\\
			&=&
			\left(\lambda_{k_1} -\lambda_{k_2}\frac{c_{k_2}}{c_{k_1}}\right)u_{k_1}\cdot g_{\gamma_{k_1}} \cdot v_{k_1}
			-S 
   		        + \sum_{k\ne k_1,k_2} \lambda_{k} u_{k}\cdot g_{\gamma_{k}}\cdot u_{k}
		\end{eqnarray*}

		We have two possibilities.  The first is that we may have succeeded in canceling all terms
		with leading monomial  $\alpha$, which contradicts the minimality of $\alpha$.  Otherwise, 
		as $\LM{S} < \alpha$,
		we have written $h$
		with no more than  $|T|-1$ terms containing $\alpha$, contradicting the
		minimality of $T$.  Thus we conclude that for such a minimal representation we must have $|T|=1$ as desired.
	\end{proof}

	\noindent Now we may proceed with the proof of our theorem.

	\begin{proof}(Theorem \ref{thm})
		Let $f$ be an element of $A$.  Then the division algorithm allows us to write:
		\begin{displaymath}
			f = r+ \tilde{f} 
		\end{displaymath}
		with $r\in O(I)$ and $\tilde(f) \in I$.   By (Lemma \ref{lem}) we see that
		we can take $r \in \tilde{O}(\G)$.  As $f$ is arbitrary in $A$, we then
		have
		\begin{displaymath}
		  A = I + \tilde{O}(\G)
		\end{displaymath}
		The theorem will follow if we can show that this sum is direct, which in turn
		will follow from showing that $r$ is unique.
		So suppose that the division algorithm produces two representations for $f$:
		\begin{displaymath}
		  f = \tilde{f} + r = \tilde{f}' + r'
		\end{displaymath}
		Then we have $\tilde{f}-\tilde{f'} \in I$ so that $r-r' \in I$.
		Now assume that $r-r'\ne 0$, then
		(Lemma \ref{lem}) shows that if there is some
		$h\in I$ such that $\LT{h} | r - r'$ then there is some $g_{\gamma} \in \G$ such 
		that $\LT{g_{\gamma}} | r- r'$.   But then, by construction of $r$ and $r'$, we know
		that there is no such $g_{\gamma}$ and we will have a contradiction by taking
		$h = r - r'$.
	\end{proof}

\begin{Proposition}\label{prop}
	Let $R$ be a commutative ring with unity, $A$ an $R$-algebra which is free as an 
	$R$-module,  $(\B,<)$ an admissible system on $A$.  Let $\alpha \in \B$ and suppose that
	$f_1, \cdots,f_n \in A$ satisfy $\LM{f_i} = \alpha$ and $\LC{f_i} \in R^{\times}$.
	Then if:
	\begin{displaymath}
		f:= \sum_{i} c_i f_i \ \ \ c_i \in R^{\ast}
	\end{displaymath}
	satisfies $\LM{f} < \alpha$ then we may write:
	\begin{displaymath}
		f = \sum_{i\ne j} d_{i,j} S_{i,j}
	\end{displaymath}
	where the $S_{i,j}$ are the $S$-polynomials about $f_i$ and $f_j$ given by:
	\begin{displaymath}
		S_{i,j} := \frac{1}{a_i} f_i - \frac{1}{a_j} f_j \ \ \ a_i := \LC{f_i},a_j:=\LC{f_j} \in R^{\times}
	\end{displaymath}
\end{Proposition}

\begin{proof}
	Because we have a cancellation of the terms of $f_i$ involving $\alpha$ we have
	that $\sum_{i} c_i = 0$.  Then:
	\begin{eqnarray*}
		\lefteqn{f= c_1 f_1 + \cdots + c_n f_n} \\
		&=& c_1 a_1 \left(\frac{1}{a_1}f_1\right) + \cdots + c_n a_n\left(\frac{1}{a_n} f_n\right) \\
		&=& c_1 a_1 \left( \frac{1}{a_1} f_1 - \frac{1}{a_2} f_2\right) + (c_1 a_1 + c_2 a_2)\left(\frac{1}{a_2}f_2 - \frac{1}{a_3} f_3\right)
			+ \cdots \\
		&& + ( c_1 a_1 + \cdots c_{n-1} a_{n-1})\left(\frac{1}{a_{n-1}}f_{n-1} - \frac{1}{a_n} f_n\right)
			+ (c_1 a_1 + \cdots c_n a_n)\frac{1}{a_n} f_n 
		\\
		& = & 
		   c_1 a_1 S_{1,2} + (c_1 a_1 + c_2 a_2) S_{2,3} + \cdots 
		\\
		& & 
			+(c_1 a_1 + \cdots + c_{n-1} a_{n-1}) S_{n-1,n} + 0 \frac{1}{a_n} f_n
	\end{eqnarray*}
	which gives the desired result.
\end{proof}

\begin{Theorem}[Buchberger]\label{buch}
	Let $R$ be a commutative ring with unity, $A$ an $R$-algebra which is free as an 
	$R$-module,  $(\B,<)$ an admissible system on $A$.  Let $I\le A$ be an ideal
	generated by a unital set:
	\begin{displaymath}
		\G:=\{ g_{\gamma} \ | \ \gamma \in \Gamma \}
	\end{displaymath}	
	for some index set $\Gamma$.   Then $\G$ is a \Grobner basis for $I$ if and only
	if all $S$-polynomials for $\G$ have zero remainder under the division algorithm.
\end{Theorem}

\begin{proof}
	We show that if all $S$-polynomials reduce to zero and  $f\in  I$ then $f$ has
	a \Grobner basis representation -- i.e a representation satisfying (Definition \ref{gbdef}).   
	As $\G$ generates $I$ we may choose a representation of $f$ as:
	\begin{equation}\label{eqnREP}
		f= \sum_{i} h_i \cdot g_{\gamma_i} \cdot h_i ' \ \ \ h_i, \ h_i'\in A \ \ \ \gamma_i \in \Gamma
	\end{equation}
	As $A$ has an $R$-basis given by $\B$ then we may write  $h_i=\sum_{k\in K} c_k \beta_k$,
	$h_i=\sum_{k'\in K'} c'_{k'} \beta'_{k'}$ with $c_k, c'_{k'} \in R^{\ast}$ and $\beta_k, \beta'_{k'} \in \B$,
	so that we have:
	\begin{displaymath}
		f= \sum_{i,k,k'} c_k c'_{k'} \beta_k \cdot g_i \cdot \beta'_{k'}
	\end{displaymath}		
	If for some representation  of $f$ as in (Equation \ref{eqnREP})
	we have for all $\beta_k\cdot g_i \cdot \beta'_{k'} \ne 0$ that
	$\LM{\beta_k\cdot g_i \cdot \beta'_{k'}} \le \LM{f}$ then we are done.
	Otherwise, let us suppose that for all such representations of $f$ we have the maximal term appearing
	$\alpha:= \max{\{\LM{h_i\cdot g_i \cdot h'_{i}}\}}$  is such that $\alpha > \LM{f}$.
	Over all such representations we may choose one so that $\alpha$ is minimal.
	We will now  produce a new representation for $f$ whose corresponding
	maximal term is strictly less than $\alpha$, thereby obtaining a contradiction.
	To this end, let us define
	$T:=\{ i | \LM{h_i\cdot g_i \cdot h'_{i}} = \alpha \}$ and:
	\begin{eqnarray*}
		g&:= &\sum_{i\in T} \LT{h_i} \cdot g_i \cdot \LT{h_i'}\\
		&=& \sum_{i\in T} \LC{h_i}\LC{h'_i} \LM{h_i} \cdot g_i \cdot \LM{h'_i}
	\end{eqnarray*}
	so that each term of  $f-g$ has leading monomial less than $\alpha$.   
	As $\G$ is assumed to be unital, we have that $a_i:=\LC{\LM{h_i}\cdot g_i \cdot\LM{h'_i}}\in R^{\times}$ so that
	we may apply to (Proposition \ref{prop}) to $g$ and write:
	\begin{equation}\label{eqnA}
		g = \sum_{i\ne j \in T} d_{i,j} S_{i,j}
	\end{equation}
	where the $S_{i,j}$ are the
	the $S$-polynomials about $\LM{h_i}\cdot g_i \cdot \LM{h'_i}$ and $\LM{h_j}\cdot g_j\cdot \LM{h'_j}$ 
	given by:
	\begin{displaymath}
		S_{i,j}:=\frac{1}{a_i}\LM{h_i}\cdot g_i \cdot \LM{h'_i} -
			 \frac{1}{a_j}\LM{h_i}\cdot g_j \cdot \LM{h'_j} 
	\end{displaymath}
	But then,  the $S_{i,j}$'s are also $S$-polynomials about $g_i$ and $g_j$, so that we have, by assumption, that
	they reduce to zero on the division algorithm, i.e. that:
	\begin{equation}\label{eqnB}
	   S_{i,j} = \sum_l \lambda_{l_{i,j}} u_{l_{i,j}}\cdot g_{l_{i,j}} \cdot v_{l_{i,j}} \ \ \ u_{l_{i,j}},v_{l_{i,j}}\in \B, \ \ 
	   \lambda_{l_{i,j}}\in R^{\ast}
	   \ \ g_{l_{i,j}} \in \G
	\end{equation}
	As $\LM{S_{i,j}} < \alpha$ we see that by
	substituting (Equation \ref{eqnB}) into (Equation \ref{eqnA})  we are able to write $g$, and thus $f$,
	in the form of (Equation \ref{eqnREP}) 
	such that
	the leading monomial of each term is strictly less than $\alpha$, our desired contradiction.
\end{proof}

\begin{Corolary}[PBW]
	  Let $\g$ be a lie algebra over $R$, a commutative ring with unity, with lie bracket $[ \ , \ ]_{\g}$.  Then:
	  \begin{displaymath}
	    \U\g = T\g/J \mbox{ where } 
			J=\left( xy-yx -[x,y]_{\g} \ | x,y \in \g \hookrightarrow T\g \right)
	  \end{displaymath}
	  is isomorphic as an $R$-module to $S\g$, the symmetric algebra on $\g$.
\end{Corolary} 
\begin{proof}
	  Choosing a well ordered basis $\{x_i \ |  \ i \in I \}$ for $\g$ then $T\g$ has a multiplicative
	  monomial basis consisting of the words of finite length in the $x_i$'s with the graded lexicographic
	  ordering.  Also $S\g$ has a basis of words of finite length written in non-decreasing order.
	  The ideal $J$ is generated by:
	  \begin{displaymath}
	    \G:= \{ g_{i,j}:=x_i x_j - x_j x_i - [x_i, x_j]_{\g}  \ | \ x_i > x_j \}
	  \end{displaymath}
	  The argument that $\G$ is a \Grobner basis is exactly
	  as in \cite{Mora} or \cite{DeGraaf} which makes use of (Theorem \ref{buch}).   
	  As the leading terms $\LT{g_{i,j}} = x_i x_j$ are monic
	  and the basis is multiplicative we have that $\G$ is a unital \Grobner basis.   
	  The corollary follows.
\end{proof}

\end{document}